\documentclass{amsart}
\usepackage{amssymb,latexsym,amsmath,amsfonts,hhline}
\usepackage{pdfsync}
\usepackage{xcolor}
\usepackage{comment,colortbl}
\usepackage{amsaddr}

\newtheorem{formula}{}[section]
\newtheorem{definition}[formula]{Definition}
\newtheorem{corollary}[formula]{Corollary}
\newtheorem{remark}[formula]{Remark}
\newtheorem{lemma}[formula]{Lemma}
\newtheorem{theorem}[formula]{Theorem}
\newtheorem{example}[formula]{Example}
\newtheorem{proposition}[formula]{Proposition}

\def\thrm{\begin{theorem}}
\def\thrml#1{\begin{theorem}\label{#1}}
\def\ethrm{\end{theorem}}
\def\prpstn{\begin{proposition}}
\def\prpstnl#1{\begin{proposition}\label{#1}}
\def\eprpstn{\end{proposition}}
\def\rmrk{\begin{remark}}
\def\rmrkl#1{\begin{remark}\label{#1}}
\def\ermrk{\end{remark}}
\def\dfntn{\begin{definition}}
\def\dfntnl#1{\begin{definition}\label{#1}}
\def\edfntn{\end{definition}}
\def\nmrt{\begin{enumerate}}
\def\enmrt{\end{enumerate}}
\def\tm#1{\item[{\rm (#1)}]}
\def\qtn{\begin{equation}}
\def\qtnl#1{\begin{equation}\label{#1}}
\def\eqtn{\end{equation}}
\def\lmm{\begin{lemma}}
\def\lmml#1{\begin{lemma}\label{#1}}
\def\elmm{\end{lemma}}
\def\crllr{\begin{corollary}}
\def\crllrl#1{\begin{corollary}\label{#1}}
\def\ecrllr{\end{corollary}}
\def\xmpl{\begin{example}}
\def\xmpll#1{\begin{example}\label{#1}}
\def\exmpl{\end{example}}

\def\css{\begin{cases}}
\def\ecss{\end{cases}}

\def\proof{\noindent{\bf Proof}.\ }

\def\oG{{\ov G}}

\def\fG{{\frak G}}

\DeclareMathOperator{\aut}{Aut}
\DeclareMathOperator{\alt}{Alt}
\DeclareMathOperator{\AGL}{AGL}

\DeclareMathOperator{\id}{id}

\DeclareMathOperator{\inv}{Inv}

\DeclareMathOperator{\orb}{Orb}

\DeclareMathOperator{\poly}{poly}

\DeclareMathOperator{\PSL}{PSL}
\DeclareMathOperator{\PGL}{PGL}
\DeclareMathOperator{\PGaL}{P\Gamma L}

\DeclareMathOperator{\soc}{Soc}

\DeclareMathOperator{\sym}{Sym}
\DeclareMathOperator{\Wr}{Wr}

\def\eprf{\hfill$\square$}

\def\grp#1{\langle {#1}\rangle}

\def\ovi#1{{\phantom{x}\hspace{-1.7mm}^{#1}}}
\def\qaq{\quad\text{and}\quad}
\def\qoq{\quad\text{or}\quad}
\def\ov{\overline}

\begin{document}

\title{Two-closure of supersolvable permutation group in polynomial time}
\author{Ilia Ponomarenko}
\address{Steklov Institute of Mathematics at  St. Petersburg, Russia\\
Mathematical Center in Akademgorodok, Novosibirsk, Russia}
\email{inp@pdmi.ras.ru}

\author{Andrey Vasil'ev}
\address{Sobolev Institute of Mathematics, Novosibirsk, Russia\\
Novosibirsk State University, Novosibirsk, Russia}
\email{vasand@math.nsc.ru}

\thanks{The work was partially supported by the RFBR grant No. 18-01-00752, and Mathematical Center in Akademgorodok}

\begin{abstract}
The $2$-closure $\ov G$ of a permutation group $G$ on $\Omega$ is defined to be the largest permutation group on $\Omega$,  having the same orbits on $\Omega\times\Omega$  as $G$. It is proved that if $G$ is  supersolvable, then
$\ov G$ can be found in polynomial time in~$|\Omega|$.
As a byproduct of our technique, it is shown that the composition factors of $\ov G$ are cyclic or alternating.
\end{abstract}

\maketitle

\section{Introduction}
It is well known that the computational problem of finding the automorphism group of a finite graph is polynomial-time equivalent to the graph isomorphism problem. None of these two problems becomes easier if the input graphs assumed to be arc-colored. One could ask whether an a priory knowledge of certain symmetries of a graph can be used to solve the former problem efficiently? In the limit case, this knowledge could include an automorphism group acting transitively on each color class of arcs. This naturally leads to the notion of the closure of a permutation group, which comes from a method of invariant relations, developed by H.~Wielandt in~\cite{Wie1969}.\medskip

In Wielandt's method, a permutation group $G\le\sym(\Omega)$  is studied by analyses of the {\em set of $m$-orbits} ($m$ is a positive integer)
$$
\orb_m(G)=\orb(G,\Omega^m),
$$
consisting of the orbits of componentwise action of~$G$ on the Cartesian power~$\Omega^m$ of the set~$\Omega$. When $m$ is sufficiently large, the group~$G$ is uniquely determined by its $m$-orbits. However, if $m$ is small, then we have to reckon with the fact that there are groups {\it $m$-equivalent} to~$G$, i.e., those having the same $m$-orbits as~$G$. For example, the only thing that can be said about $1$-equivalent groups is that they have the same orbits on~$\Omega$; one more example is that a group $2$-equivalent to a solvable group is not necessarily solvable. On the other hand, it was proved by Wielandt that the property of $G$ to be primitive ($2$-transitive, abelian, or nilpotent) is preserved with respect to the $2$-equivalence.\medskip

The {\it $m$-closure} $G\ovi{(m)}$ of a group $G\le\sym(\Omega)$ is defined to be the subgroup of $\sym(\Omega)$ leaving each $m$-orbit of~$G$ fixed,
$$
G\ovi{(m)}=\{k\in\sym(\Omega):\ s^k=s\ \, \text{for all}\ \, s\in\orb_m(G)\};
$$
the group $G$ is said to be {\it $m$-closed} if $G=G\ovi{(m)}$. The $m$-closure is $m$-equivalent to~$G$ and, in fact, does not depend on the choice of~$G$ in the class of permutation groups on~$\Omega$, that are $m$-equivalent  to~$G$. In general, $G^{(m)}=G$ for a certain $m\le |\Omega|$ and
$$
G\ovi{(1)}\ge G\ovi{(2)}\ge\cdots\ge G\ovi{(m)}=G\ovi{(m+1)}=\cdots =G.
$$
The first member $G\ovi{(1)}$ of this sequence is obviously equal to the direct
product of symmetric groups acting on the orbits of~$G$. It is far less trivial to find the $m$-closure for $m>1$. After Wielandt's pioneering work, there was some progress on the subject achieved mostly in the case of primitive groups. In \cite{LPS1988,PS1992}, the socles of the $m$-closures of primitive permutation groups were described. More recent results on invariants of the $m$-equivalence can be found in~\cite{Xu2011,Yu2019}.\medskip

The case of $2$-closure is of particular interest because of its connection with the graph isomorphism problem. Indeed, an equivalent way to define the $2$-closure $\ov G=G\ovi{(2)}$ of a group~$G\le\sym(\Omega)$  is to consider an arc-colored graph~$X$ with vertex set $\Omega$, the color classes of which are exactly the $2$-orbits of $G$;\label{270919a} this graph is nothing else than the coherent configuration associated with~$G$ in the sense of D.~Higman~\cite{Higman1975}. Then~$\ov G$ is equal to the full automorphism group $\aut(X)$ of this graph. Moreover, the natural correspondence
$$
G\to X,\quad X\to\aut(X),
$$
between the permutation groups and the corresponding graphs form a Galois correspondence; the closed objects under this correspondence are exactly the $2$-closed groups and schurian coherent configurations, see~\cite[Sec.~3.1]{EP2009} and \cite[Sec.~2.2]{CPCC}. \medskip

In this paper, we are interested in the computational complexity of {\it $2$-closure problem}: find the $2$-closure of a permutation group. As usual,  it is assumed that the input and output groups are given by sets of generators. \medskip

Modulo the recent breakthrough result by L.~Babai \cite{Bab2015}, this problem can be solved in time quasipolynomial in the degree of a given group, so the attention should be restricted to the question whether a polynomial-time algorithm exists.\medskip

From this point of view, the $2$-closure problem  was apparently first studied in~\cite{Ponomarenko1994}, where a polynomial-time algorithm  constructing the $2$-closure of a nilpotent group  was proposed. Since then the $2$-closure problem has been solved for the odd order groups~\cite{EvdokimovP2001}, for groups containing a regular cyclic subgroup~\cite{Ponomarenko2006}, and for $3/2$-transitive groups~\cite{VasilevC2018}. It should be noted that all these classes are invariant with respect to taking the $2$-closure. In the present paper, we deal with the class of supersolvable groups, in which the above property does not hold: the $2$-closure of the supersolvable group  $\AGL_1(p)$, where $p$ is a prime, is the symmetric group~$\sym(p)$. Recall that a finite group is said to be {\em supersolvable} if it contains a normal series with cyclic factors.

\thrml{200318}
The $2$-closure of a supersolvable permutation group of degree $n$ can be found in time $\poly(n)$.
\ethrm

A challenging problem is to find the $2$-closure $\ov G$ of a solvable group~$G$. The main obstacle here is that among the composition factors of~$\ov G$, simple groups other than cyclic and alternating can occur: the examples have recently been found in~\cite{Skresanov2019}. This effect does not appear for the supersolvable groups, this follows from the theorem below, obtained as a byproduct of a technique used to prove Theorem~\ref{200318}.

\thrml{200318a}
Every composition factor of the $2$-closure of a supersolvable permutation group is  a cyclic or alternating group of prime degree.
\ethrm

Theorems~\ref{200318} and~\ref{200318a} are proved in Subsections~\ref{010519d} and~\ref{040519i}, respectively. The main idea of the proof is to separate (as much as possible) the ``nonsolvable part'' of the $2$-closure $\ov G$ of a supersolvable permutation group~$G$. In this way, we define a relative closure~$K$ of~$G$ (Subsection~\ref{230919b}),  a ``solvable part'' of~$\ov G$, which is constructed with the help of the Babai-Luks algorithm~\cite{BL83}. It turns out that~$\ov G=\ov K$, and every nonsolvable section of~$\ov K$ has a very special form as a permutation group (Subsection~\ref{230919a}), which can efficiently be revealed from the corresponding section of~$G$ (Subsection~\ref{030519k}). This fact is a cornerstone of the Main Algorithm presented in Subsection~\ref{010519d}.\medskip

{\bf Notation.}

Throughout the paper, $\Omega$ is a set of cardinality~$n$.

The diagonal of $\Omega\times\Omega$ and the identity permutation of $\Omega$ are denoted by $1_\Omega$ and $\id_\Omega$, respectively.

The set of all classes of an equivalence relation $e$ on $\Omega$ is denoted by
$\Omega/e$; for a set $\Delta\subseteq\Omega$,  we put  $e_\Delta=e\cap(\Delta\times\Delta)$.

The symmetric and alternating group on $\Omega$ are denoted by $\sym(\Omega)$ and $\alt(\Omega)$, or $\sym(n)$ and $\alt(n)$ if the underlying set is fixed or irrelevant.

\section{Permutation groups: preliminaries}

Throughout the paper, all sets and groups are assumed to be finite. Our notation for permutation groups is mainly compatible with that in~\cite{DM}. In particular, given $G\le\sym(\Omega)$, we refer to the set of orbits of $G$ on $\Omega$ as $\orb(G,\Omega)$, or briefly $\orb(G)$ if $\Omega$ is fixed. The permutation group (respectively, the permutation) induced by an appropriate action of a group $G$ (respectively, of a permutation $g$) on a set $\Gamma$ is denoted by $G^\Gamma$ (respectively, $g^\Gamma$). Given $\Delta\subseteq\Omega$, we denote by $G_\Delta$ and $G_{\{\Delta\}}$ the pointwise and setwise stabilizers of~$\Delta$ in $G$, respectively. We set
$$
G^\Delta:=(G_{\{\Delta\}})^\Delta,
$$
so $G^\Delta\cong G_{\{\Delta\}}/G_\Delta$.\medskip

Given an equivalence relation~$e$ on $\Omega$, we set
$$
G_e:=\bigcap_{\Delta\in\Omega/e}G_{\{\Delta\}}.
$$
It is easily seen that the orbits of $G_e$ are subsets of the classes of~$e$. If, in addition,
$$
\orb(G_e,\Omega)=\Omega/e,
$$
then $e$ is said to be {\em normal}. \medskip

Let an equivalence relation $e$ be  $G$-invariant (in the transitive case, this means that $\Omega/e$ is an imprimitivity system for~$G$). Then $G$ acts naturally on $\Omega/e$; the permutation group induced by this action is denoted by $G^{\Omega/e}$. The kernel of the action coincides with~$G_e$; in particular, $G_e\trianglelefteq G$.  For an arbitrary  $\Delta\subseteq\Omega$, we set
$$
G^{\Delta/e}:=(G^\Delta)^{\Delta/e_\Delta}.\medskip
$$

In what follows, the $2$-closure of a permutation group~$G$ is denoted by~$\ov G$. The following statement collects some relevant properties of the group~$\ov G$; the proof is based on the Galois correspondence between permutation groups and coherent configurations.  The notation, concepts, and results used in the proof are taken from \cite{CPCC}.

 \lmml{260319a}
 Given $G\le\sym(\Omega)$ and   a $G$-invariant equivalence relation $e$ on~$\Omega$,
 \nmrt
 \tm{i} $G^{\Omega/e}$ and $\ov G\ovi{\Omega/e}$ are $2$-equivalent; in particular, $\ov G\ovi{\Omega/e}\le\ov{G^{\Omega/e}}$,
 \tm{ii} $\ov{G_e}\le\oG_e$,
 \tm{iii} if $\Delta\in\Omega/e$, then $G^\Delta$ and $\ov G\ovi\Delta$ are $2$-equivalent; in particular, $\oG\ovi\Delta\le\ov{G^\Delta}$.
 \enmrt
 \elmm
 \proof Let $\inv(G)=(\Omega,\orb_2(G))$ be the coherent configuration associated with the group~$G$. Then the $2$-closure $\ov G$ of $G$ is equal to $\aut(\inv(G))$. According to \cite[Definitions~3.1.12 and~3.1.20]{CPCC}, $\inv(G)_{\Omega/e}$ is the quotient of $\inv(G)$ modulo~$e$, and $\inv(G)_e$ is the extension of $\inv(G)$ with respect to~$e$.\medskip

(i) Applying \cite[Theorem~3.1.16]{CPCC}, we have
$$
\inv(G^{\Omega/e})=\inv(G)_{\Omega/e}=\inv(\ov G)_{\Omega/e}=\inv(\ov{G}\ovi{\Omega/e}).
$$
 
(ii) By virtue of \cite[Theorem~3.1.21]{CPCC},
$$
\ov{G_e}=\aut(\inv(G_e))\le\aut(\inv(G)_e)=(\aut(\inv(G)))_e=\ov G_e.
$$

(iii) Set $L=\ov G\ovi\Delta$. Since $\inv(G)=\inv(\oG)$, we have $\ov{G^\Delta}=\ov L$. Obviously, $\ov L\ge L$, and we are done. \eprf

\section{Permutation groups with orbits of prime cardinality}\label{030519f}

\subsection{Transitive case.}
Let $p$ be a prime and $\Delta$ a set of cardinality~$p$. Every transitive  subgroup  of $\sym(\Delta)$ contains a regular cyclic subgroup $C(\Delta)$ of order $p$. The normalizer of $C(\Delta)$  in $\sym(\Delta)$ is obviously isomorphic to $\AGL_1(p)$, and is denoted by $\AGL_1(\Delta)$.

\lmml{110117a}
Let $\Delta$ be a set of prime cardinality, $G$ a transitive subgroup of $\sym(\Delta)$, and  $C(\Delta)$ a regular subgroup of~$G$. Then either $G\le\AGL_1(\Delta)$ or~$G$ is nonsolvable and $2$-transitive. In the latter case, if $G\ge \AGL_1(\Delta)$, then $G=\sym(\Delta)$.
\elmm
\proof The first part of the lemma is well known. Indeed, as Galois proved, $G$ is solvable if and only if any two-point stabilizer of $G$ is trivial (cf., \cite[Theorem~11.6]{Wie1964}). Therefore, if $G$ is solvable and nonregular, then it is a Frobenius group, and hence $G\le\AGL_1(\Delta)$. If $G$ is nonsolvable, then it is $2$-transitive due to Burnside's theorem \cite[Theorem~3.5B]{DM}.\medskip
	
The second part can  easily be deduced from the known description of nonsolvable $2$-transitive groups with cyclic regular subgroup. Suppose that $G$ is  nonsolvable (and so $2$-transitive). Assume on the contrary that 
$$
\AGL_1(\Delta)\le G\neq\sym(\Delta).
$$ 
Let $N$  be the normalizer of $C(\Delta)$ in~$G$. It suffices to prove that the order of~$N$ is less than $|\AGL_1(\Delta)|=p(p-1)$, where $p=|\Delta|$.   This trivially holds for $G=\PSL_2(11)$, $M_{11}$, and $M_{23}$ (see, e.g.,~\cite{Atlas}). Thus we may also assume that
$$
\PSL_d(q)\le G\le\PGaL_d(q),
$$
where $d\ge2$ and  $p=(q^d-1)/(q-1)$, see \cite[Theorem~4.1]{F80}. Then $G\ge\PGL_d(q)$ and $C(\Delta)$ is conjugated to a Singer cycle of $\PGL_d(q)$ by virtue of \cite[Corollary~2]{J02}.  It follows that if $q$ is the $k$th power of a prime, then $|N|\le dkp$. On the other hand, one can verify that
$$
dk\le(q^d-1)/(q-1)-1=p-1
$$
and the equality is attained only if $d=2$ and $q=2$ or~$4$. Thus,  $|N|<p(p-1)$ unless  $p=3$ or~$5$. By the assumption, this implies that $G=\sym(3)$ or~$\sym(5)$, a contradiction. \eprf\medskip

We complete the subsection by a simple observation to be used in the proof of Theorem~\ref{300419o}.

\lmml{240419a}
Under the condition of Lemma~{\rm\ref{110117a}}, assume that $G$ is a proper subgroup of $\AGL_1(\Delta)$. Then $G$ is $2$-closed.
\elmm
\proof If $G=C(\Delta)$, then it is regular and so $2$-closed. Hence, by the hypothesis of the lemma, we may assume that $G$ is a Frobenius group. Therefore, any its irreflexive $2$-orbit is of cardinality $|G|$. It follows that no two distinct such groups have the same $2$-orbits. Thus, $G=\ov G$ as required.\eprf

\subsection{Index-distinguishable groups.}
We call a subgroup $L$ of a group $G$ {\em index-distinguishable} if $L$ is conjugated in $G$ to any subgroup of the same index. A transitive permutation group $G$ is said to be {\em distinguishable} if a point stabilizer in $G$ is an index-distinguishable subgroup of~$G$. The following lemma gives some examples of distinguishable groups. It would be interesting to get a description of primitive distinguishable groups.

\lmml{280319a}
The following permutation groups are distinguishable:
\nmrt
\tm{i} $\alt(n)$ and $\sym(n)$, $n\neq6;$
\tm{ii} a Frobenius group, in particular, any transitive subgroup of $\AGL_1(p)$ ($p$ is prime).
\enmrt
\elmm
\proof Both facts are well known. To prove (i), let $G\in\{\alt(n),\sym(n)\}$ and~$H$ a subgroup of $G$ of index $n$. Without loss of generality, we may assume that $n\ge5$. It follows that the smallest nontrivial normal subgroup of $G$ is $\alt(n)$, so the action of $G$ on the set $G/H$ of right cosets by right multiplications is faithful. Therefore, if $G=\alt(n)$ (respectively, $G=\sym(n)$), then $H$ is isomorphic to~$\alt(n-1)$ (respectively, to~$\sym(n-1)$). The rest follows from, e.g., \cite[Lemma~2.2]{W09}.\medskip

To prove~(ii), we note that any subgroup of a Frobenius group, that has the same order as a point stabilizer, meets the Frobenius kernel trivially. It follows that this subgroup is a complement of the kernel and their orders are coprime. Since all the complements are conjugated by the Schur-Zassenhaus theorem, the first statement follows. It remains to note that each transitive subgroup of $\AGL_1(p)$ is either Frobenius or regular. \eprf

\subsection{Plain groups.}\label{230919a}

Let $G\le\sym(\Omega)$. There is an equivalence relation $\sim$ on the orbits of~$G$ such that  $\Delta\sim\Gamma$ if there exists a bijection $f:\Delta\to\Gamma$ for which
\qtnl{010519h}
\alpha^{gf}=\alpha^{fg}\quad\text{for all}\ \,\alpha\in\Delta,\ g\in G.
\eqtn
i.e., the actions of $G$ on $\Delta$ and $\Gamma$ are {\em equivalent} in the sense of \cite[p.~21]{DM}. The bijection $f$ treated as a binary relation on~$\Omega$ is a $2$-orbit of~$G$. 

\dfntn
The group $G$ is said to be plain  if given $\Delta,\Gamma\in\orb(G)$,
$$
\Delta\times\Gamma\in\orb_2(G)\qoq \Delta\sim\Gamma.
$$
\edfntn

Clearly, the two above conditions for $\Delta$ and $\Gamma$ are satisfied simultaneously if and only if $\Delta$ and $\Gamma$ are both singletons; in particular, an identity permutation group is always plain.\medskip

In the following statement, we establish a sufficient condition for a permutation group to be plain. Recall that a permutation group is said to be $\frac{1}{2}$-{\em{transitive}} if all its orbits are of the same size. A permutation group is said to be {\em quasiprimitive} if every nontrivial normal subgroup of it is transitive. It is easily seen that primitive groups  and transitive simple  groups are quasiprimitive.

\lmml{010418a}
A $\frac{1}{2}$-transitive permutation group $G$ is plain if every transitive constituent of~$G$ is quasiprimitive and distinguishable.
\elmm

\proof Let $\Delta$ and $\Gamma$ be two distinct $G$-orbits. We may assume that both of them are non-singletons, for otherwise  $\Delta\times\Gamma\in\orb_2(G)$. Recall that $G^{\Delta\cup\Gamma}\cong G/(G_\Delta\cap G_\Gamma)$ and set
$$
f_\Delta: G^{\Delta\cup\Gamma}\rightarrow G^\Delta\cong G/G_\Delta\qaq f_\Gamma: G^{\Delta\cup\Gamma}\rightarrow G^\Gamma\cong G/G_\Gamma
$$
to be  the natural epimorphisms from $G^{\Delta\cup\Gamma}$ onto the transitive constituents $G^\Delta$ and $G^\Gamma$, respectively.\medskip

First, we assume that one of these epimorphisms, say $f_{\Delta}$, is not injective. Then $$
G_\Delta=\ker(f_\Delta)>G_\Delta\cap G_\Gamma.
$$
Therefore, $G^\Gamma\cong G/G_\Gamma$ includes the nontrivial normal subgroup
$$
(G_\Delta)^\Gamma\cong G_\Delta G_\Gamma/G_\Gamma\cong G_\Delta/(G_\Delta\cap G_\Gamma).
$$
Since $G^\Gamma$ is quasiprimitive,  $(G_\Delta)^\Gamma$ is transitive. It follows that $\Delta\times\Gamma\in\orb_2(G)$.\medskip

Thus, we may assume that both $f_\Delta$ and $f_\Gamma$ are isomorphisms. This implies that the composition $f= f_\Delta^{-1}f^{}_\Gamma$ is a group isomorphism from $G^\Delta$ onto~$G^\Gamma$. Since $G^\Gamma$ is distinguishable, the image $H^f$ of a point stabilizer $H$ in $G^\Delta$ is conjugated to some point stabilizer in $G^\Gamma$. It follows that $H^f$ is a point stabilizer itself, so the actions of $G$ on $\Delta$ and $\Gamma$ are equivalent by \cite[Lemma~1.6B]{DM}, and we are done. \eprf

\crllrl{290419a}
Suppose that  all orbits of a permutation group $G$ are of the same prime cardinality, and all nonsolvable transitive constituents of $G$ are alternating or symmetric. Then $G$ is plain.
\ecrllr
\proof In view of the previous lemma it suffices to check that the transitive constituents of $G$ are distinguishable. However, this immediately follows from Lemmas~\ref{110117a} and~\ref{280319a}.\eprf\medskip

When $G$ is a plain group, we say that $\sim$ and the~$f$ are  the {\em standard equivalence relation} and {\it plain bijections}, respectively. In what follows, given a class $\Lambda$ of $\sim$, the union of all orbits of $G$ contained in $\Lambda$ is denoted by $\Omega_\Lambda$.

\lmml{240318a}
Let $G\le\sym(\Omega)$ be a plain group. Then
\nmrt
\tm{i} for each $\Lambda\in\Omega/\sim$  and  each~$\Delta\in\Lambda$, the epimorphism $G^{\Omega_\Lambda}\to G^\Delta$ is injective,
\tm{ii} if $G^\Delta$ is primitive and nonregular for all~$\Delta\in\orb(G)$, then there is at most one plain bijection $f:\Delta\to\Gamma$ for fixed  $\Delta,\Gamma\in\orb(G)$,
\tm{iii} if $G^\Delta$ is a distinguishable nonabelian simple group for all ~$\Delta\in\orb(G)$, then $G=\prod_{\Lambda\in\Omega/\sim}G^{\Omega_\Lambda}$.
\enmrt
\elmm

\proof To prove (i), let $g\in G$ satisfy $g^\Delta=\id_\Delta$. We should verify that $g^{\Omega_\Lambda}=\id_{\Omega_\Lambda}$. To this end, let $\beta\in\Omega_\Lambda$. Then there exists $\alpha\in\Delta$ such that $\beta=\alpha^f$ for a suitable plain bijection~$f$. It follows that
$$
\beta^g=(\alpha^f)^g=(\alpha^g)^f=\alpha^f=\beta,
$$
as required.\medskip

To prove (ii), assume on the contrary that $f':\Delta\to\Gamma$ is a plain bijection different from~$f$. Then $f'f^{-1}\ne 1_\Delta$ is an $G^\Delta$-invariant binary relation of cardinality~$\Delta$. It follows that there exists $\alpha\in\Delta$ such that $\beta:=\alpha^{f'f^{-1}}$ is different from~$\alpha$ and
$$
(G^\Delta)_\alpha=((G^\Delta)_\alpha)^{f'f^{-1}}=(G^\Delta)_{\alpha^{f'f^{-1}}}=(G^\Delta)_\beta.
$$
Since the group $G^\Delta$ is primitive, this is possible only if it is regular \cite[Proposition~8.6]{Wie1964}, a contradiction.\medskip

Let us prove (iii). By statement~(i), for each $\Lambda\in\Omega/\sim$  and  each~$\Delta\in\Lambda$, the group $G^{\Omega_\Lambda}$ is isomorphic to~$ G^\Delta$ and hence is nonabelian simple. Therefore it suffices to verify that given $\Delta,\Gamma\in\orb(G)$ such that $\Delta\not\sim\Gamma$,
$$
G^{\Delta\cup\Gamma}=G^\Delta\times G^\Gamma.
$$
Assume the contrary. Then the simplicity of $G^\Delta$ and~$G^\Gamma$ implies that the restriction homomorphisms $G^{\Delta\cup\Gamma}\to G^\Delta$ and $G^{\Delta\cup\Gamma}\to G^\Gamma$ are isomorphisms.  Since $G^\Delta$ and~$G^\Gamma$ are distinguishable, they are equivalent (see the proof of Lemma~\ref{010418a}). It follows that $\Delta\sim\Gamma$, a contradiction. \eprf

\subsection{Auxiliary lemma}

The following lemma establishes an important property of normal plain subgroups of a transitive group; it will be used in the proofs of Theorems~\ref{300419o} and~\ref{230419b}.

\lmml{010418b}
Let  $e$ be an equivalence relation on a set $\Omega$. Assume that a transitive group  $G\le\sym(\Omega)$ and its normal subgroup $L$ satisfy the following conditions:
\nmrt
\tm{i} $e$ is $G$-invariant and $\orb(L)=\Omega/e$,
\tm{ii} $L$ is plain,
\tm{iii} each transitive constituent of~$L$ is primitive and nonregular.
\enmrt Then  the standard equivalence relation for~$L$ does not depend (for a fixed~$e$) on $G$ and $L$, and depends only on the $2$-orbits of~$G$.
\elmm

\proof It suffices to verify that given $\Delta,\Gamma\in\orb(L)$ and any $r\in\orb_2(G)$ intersecting $\Delta\times\Gamma$,
$$
\Delta\times\Gamma\in\orb_2(L)\quad\Leftrightarrow\quad \Delta\times\Gamma\subseteq r.
$$
The implication $\Rightarrow$ follows from the fact that each $2$-orbit of $L$ is contained in a unique $2$-orbit of~$G$. Conversely, suppose on the contrary that $\Delta\times\Gamma$ is contained in a certain $r\in\orb_2(G)$ but is not a $2$-orbit of~$L$. Since $L$ is plain, this implies that there is a plain bijection
$$
f:\Delta\to\Gamma.
$$
Taking into account that $L\trianglelefteq G$, we see that for each  $g\in G$, the bijection
$$
g^{-1}fg:\Delta^g\to\Gamma^g
$$
is also plain. Moreover, the hypothesis on the transitive constituents of~$L$ and Lemma~\ref{240318a}(ii) imply that either $f=g^{-1}fg$, or $(\Delta\times\Gamma)\cap g^{-1}fg=\varnothing$. This shows that if~$r'$ is the union of all $g^{-1}fg$, then
$$
r'\cap (\Delta\times\Gamma)=\{f\}\ne\Delta\times\Gamma,
$$
unless $|\Delta|=|\Gamma|=1$, i.e., $e=1_\Omega$. But the latter is impossible
in view of~(i) and~(iii).\medskip

On the other hand, $r'$ is obviously a $G$-invariant relation intersecting $r$. Since $r$ is a $2$-orbit of $G$, this implies that $r\subseteq r'$. But then
$$
\Delta\times\Gamma\subseteq r\subseteq r',
$$
a contradiction.\eprf

\section{Permutation groups: flags, majorants, and relative closures}\label{210419b}

Throughout this section, $G\le\sym(\Omega)$ and $F=\{e_i:\ i=0,\ldots,m\}$ is a family  of $G$-invariant equivalence relations.

\subsection{Flags}
 We say that $F$  is a {\it $G$-flag} of  {\it length} $m$ if 
 \qtnl{170419a}
1_\Omega=e_0\subset e_1\subset\cdots\subset e_m\qaq \Omega/e_m=\orb(G,\Omega).
 \eqtn
 In this case, each class of the equivalence relation $e_i$  is contained in a uniquely determined $G$-orbit and forms a block in the corresponding transitive constituent.\medskip

The flag $F$ is said to be  {\em normal} if the equivalence relation $e_i$ is normal  with respect to~$G$ for all~$i$. In this case,  the equivalence relation on $\Omega/e_{i-1}$ with classes  $\Delta/e_{i-1}$, $\Delta\in\Omega/e_i$, is  normal and $G^{\Omega/e_{i-1}}$-invariant.  A normal $G$-flag $F'$ {\it extends} $F$ if $F'\supseteq F$; when $F'\ne F$, we say that $F'$ {\it strictly extends}~$F$.  A normal flag $F$ is said to be {\it maximal} if no normal $G$-flag strictly extends $F$. Thus $F$ is maximal  if  there is no normal equivalence relation lying strictly between $e_{i-1}$ and~$e_i$ for all~$i\ge 1$.\medskip

Let 
\qtnl{290319a}
1=L_0<L_1<\cdots<L_m=G
\eqtn
be  a normal series for the group $G$. Given $i\in\{0,1,\ldots,m\}$, denote by $e_i$  the equivalence relation on $\Omega$ such that
$$
\Omega/e_i=\orb(L_i,\Omega),
$$
in particular, the classes of~$e_m$ are the orbits of~$G$. The family $F=\{e_i\}_{i=0}^m$ is obviously a normal $G$-flag. In addition, if $(L_i)^{\Omega/e_{i-1}}$ is a minimal normal subgroup of $G^{\Omega/e_{i-1}}$,  $i=1,\ldots,m$, then the flag $F$ is maximal. Thus, the following statement holds.

\lmml{121219a}
For every permutation group $G$ there exists a maximal normal $G$-flag.
\elmm

Given a $G$-invariant set $\Delta$, one can define a $G^\Delta$-flag
$$
F_\Delta=\{e_\Delta:\ e\in F\}
$$
and given an equivalence relation  $e\in F$, one can define  a $G^{\Omega/e}$-flag
$$
F_{\Omega/e}=\{e'_{\Omega/e}:\ e'\in F,\ e'\supseteq e\},
$$
where $e'_{\Omega/e}$ is the equivalence relation the classes of which are $\Delta/e$ with $\Delta\in\Omega/e'$. Of course, in both cases, the length of the new flag can be smaller than~$m$. One can see that the flags $F_\Delta$ and $F_{\Omega/e}$ are  normal  if so is~$F$.\medskip

The following statement directly follows from the definitions.

\lmml{170419b}
Let $G$ and $K$  be permutation groups on the same set, and let $F$ be a $G$-flag. Then
\nmrt
\tm{i} if $G$ and $K$ are $2$-equivalent, then $F$ is a $K$-flag;
\tm{ii} if $K\ge G$ and $F$ is a $K$-flag, then  $F$ is normal as $K$-flag whenever $F$ is a normal as $G$-flag.
\enmrt
\elmm

\subsection{Sections}
Let $F$ be an arbitrary  $G$-flag of length~$m$.  Given $i\in\{1,\ldots,m\}$, the group $G_{e_i}$ acts on the set $\Omega/e_{i-1}$. Denote by  $\Omega_i$ the union of all non-singleton orbits of the induced permutation group. This set is definitely nonempty if the equivalence relation $e_i$ is  normal. In the latter case, the  group
$$
S:=(G_{e_i})^{\Omega_i}.
$$
is called a {\it section} of $G$ with respect to~$F$, or briefly, an $F$-section of~$G$. We put $i(S)=i$, $e_S:=e_i$,  and $\Omega_S:=\Omega_i$.

\lmml{170419c}
If $S$ is a section of $G$ with respect to a maximal normal flag~$F$, then  the action of $S$ on any nonempty $G^{\Omega_S}$-invariant set is faithful. In particular, the transitive constituents of $S$ are isomorphic groups.
\elmm
\proof Let $\Delta\subseteq \Omega_S$ be a nonempty $G\ovi{\Omega_S}$-invariant set. Clearly, we may suppose that $\Delta\neq\Omega_S$. Denote by $\ov e$ the equivalence relation on $\Omega_S$ such that
$$
\Omega_S/\ov e=\orb(S_\Delta,\Omega_S),
$$
where $S_\Delta$ is the pointwise stabilizer of $\Delta$ in $S$.  From the choice of $\Delta$, it follows that $S_\Delta$ is normal in $G^{\Omega_S}$, so $\ov e$ is $G^{\Omega_S}$-invariant. Assume on the contrary that the action of $S$ on $\Delta$ is not faithful. Then
\qtnl{170419d}
1_{\Omega_S}\subset \ov e\subset e_S.
\eqtn
Now let $i=i(S)$. Denote by $e$ the equivalence relation on $\Omega$ such that the restriction of $e$ to $\Omega/{e_{i-1}}$ coincides with $\ov e$ inside $\Omega_i$ and has singleton classes outside $\Omega_i$. Then~$e$ is $G$-invariant and normal. Moreover,  in view of~\eqref{170419d},
$$
e_{i-1}\subset e\subset e_i
$$
in contrast to the maximality of the flag~$F$.\eprf\medskip

From Lemma~\ref{170419c}, it follows that if $S$ is a nonsolvable section with respect to a maximal normal flag then all transitive constituents of~$S$ are nonsolvable.

\subsection{Majorants}\label{230919b}
Let $F$ be an arbitrary  $G$-flag of length~$m$.  Assume first that the group $G$ is transitive. For each $i$, choose a class $\Delta_i$ of the equivalence relation~$e_i$ so that
$$
\Delta_0\subset\Delta_1\subset\cdots\subset\Delta_m.
$$
In particular, $\Delta_0$ is a singleton and $\Delta_m=\Omega$. Set $\ov\Delta_i=\Delta_i/e_{i-1}$, $i=1,\ldots,m$, and identify $\Omega$  with the Cartesian product  of the~$\ov\Delta_i$,
\qtnl{270419i}
\Omega={\ov\Delta_1}\times\cdots \times {\ov\Delta_m}.
\eqtn
Note that this identification is not canonical and depends on the choice of suitable permutations belonging to~$G$. \medskip

Under identification~\eqref{270419i}, the group $G$ can naturally be treated as a subgroup of the iterated (imprimitive) wreath product
\qtnl{27042019j}
\Wr(G;\Delta_0,\ldots,\Delta_m)=G^{\ov\Delta_1}\wr G^{\ov\Delta_2}\wr\cdots\wr G^{\ov\Delta_m},
\eqtn
where for $i=1,\ldots, m$, we set  $G^{\ov\Delta_i}$ to be the permutation group induced by the action of $G^{\Delta_i}$ on $\ov\Delta_i$. In the spirit of~\cite{EP2003}, the group~\eqref{27042019j}  is called an {\it $F$-majorant} for~$G$ and  is denoted by $W_F(G)$. It should be stressed that a different  choice of the permutations used for the above identification can lead to different majorants. But once the identification has been fixed, the $F$-majorant is well defined. (cf. \cite[Lemma~3.4]{EP2003})\medskip

In the general case, when the group $G$ is not necessarily transitive, the majorant is defined to be the direct product of the majorants for the transitive constituents of~$G$,
$$
W_F(G)=\prod_{\Delta\in\orb(G)} W_{F_\Delta}(G^\Delta),
$$
where the direct product acts on the disjoint union of the underlying sets of the factors. As is easily seen from the definition, $G\le W_F(G)$.

\dfntn
The group $K=\ov G\cap W_F(G)$ is called a relative closure of the group~$G$ with respect to the flag $F$ $($and the majorant $W_F(G))$.
\edfntn

We conclude this subsection collecting some elementary properties of a majorant and relative closure of $G$.

\lmml{251019a}
Let $G\le\sym(\Omega)$, $F$ a $G$-flag, $W=W_F(G)$ a majorant for $G$, and~$K$ the relative closure of $G$ with respect to~$F$. Then
\nmrt
\tm{i} $\ov G=\ov K$,
\tm{ii} $\ov{K}\cap W_F(K)=K$,
\tm{iii} if $G$ is solvable, then $W$ and $K$ are solvable.
\enmrt
\elmm
\proof We have $G\le K\le\ov G$ and hence $K$ is  $2$-equivalent to~$G$.  This yields~(i). \medskip

Under identification~\eqref{270419i}, every orbit $\Delta$ of $G$ is represented as the direct product $\Delta={\ov\Delta_1}\times\cdots \times {\ov\Delta_m}.$ By the definition of $F$-majorant,
\qtnl{031219z}
G\ovi{\ov\Delta_i}=W\ovi{\ov\Delta_i},\quad i=1,\dots,m.
\eqtn
This immediately implies (iii).\medskip

Next,  $G\le K\le W$ yields $W\le W_F(K)\le W_F(W)$. In view of equalities~\eqref{031219z},  $W=W_F(W)$. Thus, $\ov{K}\cap W_F(K)=\ov{G}\cap W=K$, and (ii) holds. \eprf

\subsection{Standard representation}\label{131219o}
In what follows, we use a {\it standard representation} for the permutations belonging to the group $W:=W_F(G)$, where the $G$-flag $F$ is assumed to be normal. Namely, such a representation of a permutation $k\in W$ is given by a family $\{k_\Delta\}$, where $k_\Delta\in G^\Delta$ and $\Delta$ runs over the elements of the sets
$$
O_i=\orb(G_{e_i},\Omega/e_{i-1}),\qquad  i=1,\ldots,m;
$$
 the permutation~$k_\Delta$ is called a {\it $\Delta$-coordinate} of~$k$ and is defined as follows.\medskip

Let $i\ge 1$ and $\Delta\in O_i$. Denote by $\Lambda$ the orbit of $G^{\Omega/e_{i-1}}$ that contains $\Delta$. Then
$$
k^\Lambda\in W_{F'}(G')\le (G')^\Delta\wr G^{\Lambda/e_i},
$$
where $F'$ and $G'$ are the restrictions of $F_{\Omega/e_{i-1}}$ and $G^{\Omega/e_{i-1}}$ to~$\Lambda$, respectively. Now, the $\Delta$-coordinates of the permutation~$k$ are obtained from the representation of $k$ as an element of the wreath product on the right-hand side of the above inclusion:
$$
k^\Lambda=(\{k_\Delta:\ \Delta\in O_\Lambda\}; k^{\Lambda/e_i}).
$$
where $O_\Lambda=\{\Delta\in O_i:\ \Delta\subseteq\Lambda\}$.\medskip

We note that by an inductive argument, the permutation $k^{\Lambda/e_i}$ can also be written in terms of the $\Delta$-coordinates of~$k$, corresponding to the indices~$i,\ldots,m$. Since $\Omega/e_{i-1}$  is a disjoint union of $\Lambda\in\orb(G^{\Omega/e_{i-1}})$,
$$
W^{\Omega/e_{i-1}}=\prod_{\Lambda\in\orb(G^{\Omega/e_{i-1}})}W^\Lambda.
$$
It follows that the action of $k$ on $\Omega/e_{i-1}$ can be written in the form $(\{k_\Delta\}; k^{\Omega/e_i})$, where this time $\Delta$ runs over the whole set~$O_i$.  The multiplication of the permutations given in this form is performed as in the case of wreath product. In particular,
\qtnl{040519a}
(\{l_\Delta\}; \id_{\Omega/e_i})\cdot (\{k_\Delta\}; k^{\Omega/e_i})=(\{l_\Delta k_\Delta\}; k^{\Omega/e_i}).
\eqtn

The $\Delta$-coordinates taken for a plain section satisfy some additional relations. Namely, the following statement is true.

\lmml{040519b}
Let $F$ be a normal $G$-flag. Then 
for every  plain $F$-section~$S$ of $G$, 
$$
\Delta\sim\Gamma\quad\Rightarrow\quad k_\Gamma=f^{-1}k_\Delta f,
$$
for all $k\in G_{e_S}$ and $\Delta,\Gamma\in\orb(S)$,  where $\sim$ is  the standard equivalence relation for~$S$ and $f:\Delta\to\Gamma$ is a plain bijection. 
\elmm
\proof  Let $k\in G_{e_S}$. Then  in the standard representation, $k=(\{k_\Delta\}; 1)$, where $1$ denotes $\id_{\Omega/e_S}$. By the definition of plain bijection, this implies that for all $\beta\in\Gamma$,
$$
\beta^{k_\Gamma}=\beta^{(\{k_\Gamma\};1)}=
(\beta^{f^{-1}})^{f(\{k_\Gamma\};1)}=
(\beta^{f^{-1}})^{(\{k_\Gamma\};1)f}=
(\beta^{f^{-1}})^{k_\Delta f},
$$
as required.\eprf

\section{The relative closure  of a supersolvable permutation group}\label{040519u}

Throughout this section, $K\le\sym(\Omega)$ is the relative closure of a supersolvable group~$G$, and  $\ov K=K^{(2)}$ is the $2$-closure of~$K$. It should be noted that the group~$K$ is not necessarily supersolvable.

\subsection{The orbits of a section of~$K$.}
Let $F$ be a maximal normal $G$-flag of length~$m\ge 1$. Observe that $F$ is also a  normal $K$- and $\ov K$-flag, see Lemma~\ref{251019a}.  

\thrml{181219a}
Let $S$ be an $F$-section of $K$. Then there exists a prime $p$ such that $|\Delta|=p$ for all $\Delta\in\orb(S)$.
\ethrm
\proof Let $T$ be the $F$-section of~$G$ with $i(T)=i(S)=:i$. The normality of~$F$ as a $G$-flag and as a $K$-flag yields
$$
\orb(T)=\Omega_T/e=\Omega_S/e=\orb(S),
$$
where $e=(e_i)_{\Omega/e_{i-1}}$. Thus it suffices to show that $|\Delta|=p$ for all $\Delta\in\orb(T)$.\medskip

Let $1<N\le T$ be a  minimal normal subgroup of $G^{\Omega_T}$. Then the orbits of~$N$ define a $G^{\Omega_T}$-invariant equivalence relation. Therefore, the maximality of~$F$ as a $G$-flag implies that
$$
\orb(N)=\orb(T).
$$ 
Moreover, since $G$ and so~$G^{\Omega_T}$ is supersolvable, $N$ is a cyclic group of prime order~$p$. Thus any orbit of $N$ and hence of~$T$ is of cardinality~$p$, as required.\eprf

\crllrl{181219z}
Let $F'$ be a normal $K$-flag extending the flag~$F$. Then 
$$
F'_\Lambda=F^{}_\Lambda,\qquad\Lambda\in\orb(K).
$$
\ecrllr
\proof  Let $\Lambda\in\orb(K)$. From Theorem~\ref{181219a}, it follows that given $\Delta\in\Lambda/e_i$ and $\Delta'\in\Delta/e_{i-1}$ ($1\le i\le m$),
$$
|\Delta|=|\Delta'|\qoq |\Delta|=p|\Delta'|.
$$
On the other hand, any two classes of every  $K^\Lambda$-invariant equivalence relation, being blocks of $K^\Lambda$, have the same cardinality. Thus the normal flag $F_\Lambda$ is maximal. Since~$F'_\Lambda$ is a normal flag extending $F_\Lambda$, we are done.\eprf\medskip

From Corollary~\ref{181219z}, it follows that if $G$ is transitive, then every maximal normal $G$-flag is also maximal normal $K$-flag. In the intransitive case, we cannot guarantee the maximality.  The following statement enables us to control it (cf. Lemma~\ref{170419c}).

\lmml{181219b}
Let $F'$ be a nonmaximal normal $K$-flag extending the flag~$F$. Then there is an equivalence relation $e=e'_{i-1}$ of~$F'$, a non-singleton  $K\ovi{\Omega/e}$-orbit $\Lambda$, and a minimal normal subgroup $N$ of~$ K\ovi{\Omega/e}$, contained in $(K_{\ov{e}})\ovi{\Omega/e}$ with $\ov e=e'_i$, such that 
$$
N^\Lambda=\{\id_\Lambda\}\qaq N\ovi{\Lambda'}\ne\{\id_{\Lambda'}\},
$$ 
where $\Lambda'$ is the complement of $\Lambda$ in $\Omega/e$.
\elmm
\proof Let $F'=\{e'_i:\ i=0,\ldots,m'\}$, where $m'\ge m$. By the assumption, there exists an integer $1\le i\le m'$ and an equivalence relation $e^*$ on $\Omega$ such that 
$$
\{e'_0,\ldots,e'_{i-1},e^*,e'_i,\ldots,e'_{m'}\}
$$ 
is a normal $K$-flag extending~$F'$. Setting $e=e'_{i-1}$ and $\ov e=e'_i$, we have  $e\subsetneq e^*\subsetneq \ov e$. This implies that  
$$
1_\Lambda=e^{}_\Lambda\subseteq e^*_\Lambda\subseteq \ov e_\Lambda,\qquad
\Lambda\in\orb(K\ovi{\Omega/e}),
$$
and for at least one non-singleton orbit $\Lambda$, the right-hand side inclusion is strict. On the other hand, the maximality of $F$ implies that the normal flag $F_\Lambda$ is maximal. Thus by Corollary~\ref{181219z}, we conclude that 
$$
e^*_\Lambda=1^{}_\Lambda.
$$
Since the flag $F'$ is normal, the orbits of $(K_{e^*})\ovi\Lambda$ are singletons. Moreover, the group $(K_{e^*})\ovi{\Lambda'}$ is nontrivial, for otherwise $e^*=e$. This proves the required statement for any minimal normal subgroup $N$ of $K\ovi{\Omega/e}$, contained in $(K_{e^*})\ovi{\Omega/e}\le(K_{\ov e})\ovi{\Omega/e}$.\eprf

\crllrl{181219i}
For every maximal normal $G$-flag  $F$, there exists a maximal normal $K$-flag extending~$F$.
\ecrllr

\subsection{The constituens of a section of $K$.} From now on, we always assume that~$F$ is a maximal normal $K$-flag. Let $S$ be an $F$-section of $K$. Since $F$ is a normal $\ov{K}$-flag, the $F$-section $\ov S$ of $\ov K$ satisfying $\Omega_{\ov S}=\Omega_{\vphantom{\ov S}S}$ is well defined. Clearly, 
$$
\orb(\ov S)=\orb(S)\qaq \ov S\le S^{(2)}.
$$
It should be noted that $\ov S$ and $S^{(2)}$ do not necessarily coincide  (cf. Lemma~\ref{260319a}).

\thrml{300419o}
Let $S$ be an $F$-section of $K$. Then
\nmrt
\tm{i}  given $\Delta\in\orb(S)$,
\nmrt
\tm{a} if $\ov S$ is solvable, then $\ov K\ovi{\Delta}=K^\Delta\le\AGL_1(\Delta)$,
\tm{b}  if $\ov S$ is nonsolvable, then $|\Delta|\ge 5$, $K^\Delta=\AGL_1(\Delta)$,  $\ov K\ovi{\Delta}=\sym(\Delta)$, and $\soc(\ov S\ovi{\Delta})=\alt(\Delta)$,
\enmrt
\tm{ii} $S$, $\soc(\ov S)$, and $\ov S$ are plain groups; if $\ov S$ is nonsolvable, then the corresponding standard equivalence relations of $\soc(\ov S)$ and $\ov S$ coincide.
\enmrt
\ethrm
\proof  To prove  (i), let $\Delta\in\orb(S)$. By statements (i) and (iii) of Lemma~\ref{260319a}, the groups $\ov K\ovi{\Delta}$ and $K^\Delta$ are $2$-equivalent.  Moreover, both of them are of prime degree~$p$ by Theorem~\ref{181219a}.  Now if $\ov S$ is nonsolvable, then $p\ge 5$ and $\ov K\ovi{\Delta}\ge \ov S\ovi{\Delta}$ is also nonsolvable (Lemma~\ref{170419c}).  By Lemma~\ref{110117a}, this implies that $\ov K\ovi{\Delta}$ and hence $K^\Delta$ is $2$-transitive; in particular, $K^\Delta=\AGL_1(\Delta)$. Taking into account that $\ov K\ovi{\Delta}\ge K^\Delta$, we  conclude by the same lemma that $\ov K\ovi{\Delta}=\sym(\Delta)$.\medskip

Let $\ov S$ be solvable. If, in addition, $\ov K\ovi{\Delta}$  is solvable, then the required statement follows from Lemmas~\ref{110117a} and~\ref{240419a}. Assume that  $\ov K\ovi{\Delta}$  is nonsolvable. Then $p\ge 5$, and each normal transitive subgroup of  $\ov K\ovi{\Delta}$ coincides with $\alt(\Delta)$ or $\sym(\Delta)$ (Lemma~\ref{110117a}). Therefore, the group $\ov S\ovi{\Delta}=(\ov K_e)^\Delta$ being a nontrivial normal subgroup of $\ov K\ovi{\Delta}$ contains $\alt(\Delta)$ as a subgroup. Then $\ov S$ is not solvable in contrast to the hypothesis. This completes the proof of~(i).\medskip

The first part of statement~(ii) follows from Theorem~\ref{181219a}, statement~(i), and Corollary~\ref{290419a}. To prove the second part, we apply Lemma~\ref{010418b}, setting
$$
\Omega=\Omega_S,\quad e=e_S,\quad G=\ov K\ovi{\Omega_S},\quad L\in\{\soc(\ov S),\ov S\}
$$
in the hypothesis of this lemma. Obviously, $e_S$ is $\ov K\ovi{\Omega_S}$-invariant. Since $\ov{S}$ is nonsolvable, we have $\soc(\ov S\ovi{\Delta})=\alt(\Delta)$ for each $\Delta\in\orb(S)$ (see statement (i)(b) above). It follows that condition~(iii) of Lemma~\ref{010418b} holds. Moreover, by Theorem~\ref{181219a}, 
$$
\orb(\soc(\ov S))=\orb(S)=\orb(\ov{S}),
$$ 
which proves  condition (i) of that lemma. Finally, condition~(ii) follows from the first part of statement~(ii).  Now Lemma~\ref{010418b} yields that the standard equivalence relations for $\soc(\ov S)$ and $\ov{S}$ are the same, as required.\eprf

\subsection{A reduction lemma.}\label{040519i} A key point in our arguments is the following statement, which enables us to lift certain permutations from sections of $\ov K$.

\lmml{040519c}
Let $k\in \ov K$, and $\ov S$ be an $F$-section of~$\ov K$. Assume that $k_\Delta\not\in K^\Delta$ for some $\Delta\in\orb(\ov S)$. Then $\ov S$ is nonsolvable and there exists $l\in\ov K$ such that $l^{\Omega_S}\in\soc(\ov S)$, and
\qtnl{060519a}
(lk)^{\Omega/e_S}=k^{\Omega/e_S}\qaq   (lk)_\Delta\in K^\Delta\ \,\text{for all}\ \,\Delta\in\orb(\ov S).
\eqtn
\elmm
\proof The section $\ov S$ is nonsolvable due to statement~(i)(a) of Theorem~\ref{300419o}. By statement~(i)(b) of the same theorem, this yields 
$$
\soc(\ov S\ovi{\Delta})=\alt(\Delta)\qaq K^\Delta=\AGL_1(\Delta)
$$ 
for all $\Delta\in\orb(\ov S)$. Since $\sym(\Delta)=\langle\alt(\Delta),\AGL_1(\Delta)\rangle$, it follows that there exist $u(\Delta)\in\soc(\ov S\ovi{\Delta})$ and $v(\Delta)\in K^\Delta$ such that
$$
k_\Delta=u(\Delta)v(\Delta).
$$
The group $\ov S$ is plain (Theorem~\ref{300419o}(ii)) and the transitive constituents of $\ov S$ are primitive and nonregular (Theorem~\ref{300419o}(i)(b)). Hence for all  $\Delta,\Gamma\in\orb(\ov S)$
$$
\Delta\sim\Gamma\quad\Rightarrow\quad f^{-1}u(\Delta)v(\Delta)f=u(\Gamma)v(\Gamma),
$$
where $\sim$ is the standard equivalence relation for~$\ov S$ and $f:\Delta\to\Gamma$ is a plain bijection (Lemma~\ref{040519b}). \medskip

In each class of the standard equivalence relation, fix an $\ov S$-orbit $\Delta_0$. For  every $\Delta\sim\Delta_0$,  take  $u(\Delta)$ and $v(\Delta)$ so that $u(\Delta)=f^{-1}u(\Delta_0)f$. Theorem~\ref{300419o}(ii) implies that the standard equivalence relations for $\ov S$ and $\soc(\ov S)$ are the same, so applying Lemma~\ref{240318a}(iii), we obtain
$$
\prod_{\Delta\in\orb(\ov S)}u(\Delta)^{-1}\in\soc(\ov S),
$$
where the factor $u(\Delta)^{-1}$ is interpreted as the permutation on $\Omega_S$ with support~$\Delta$.
It follows that $\ov K$ contains  an element $l$ whose $\Delta$-coordinates are equal to $u(\Delta)^{-1}$ for all $\Delta\in\orb(\ov S)$, and also $l^{\Omega/e_S}=\id_{\Omega/e_S}$. In particular, $l^{\Omega_S}\in\soc(\ov S)$ and the left-hand side equality in~\eqref{060519a} holds. Furthermore, using formula~\eqref{040519a}, we have
$$
(lk)_\Delta=l_\Delta k_\Delta=u(\Delta)^{-1}(u(\Delta)v(\Delta))=v(\Delta)\in K^\Delta,
$$
as required.\eprf

\crllrl{060519c}
Let  $\ov S$ be an $F$-section of~$\ov K$. Then every  $\ov k\in\ov S$ can be lifted to $k\in\ov K$ such that $k_\Delta\in K^\Delta$ for all $\Delta\not\in\orb(\ov T)$, where $\ov T$ runs over $F$-sections of~$\ov K$ with $i(\ov T)\ne i(\ov S)$.
\ecrllr
\proof We have $\ov k=k^{\Omega_S}$ for some $k\in\ov K_{e_S}$. Let us apply Lemma~\ref{040519c} consecutively to each section $\ov T$ of $\ov K$, $i(\ov T)=i(\ov S)-1,\ldots,1$, for which the hypothesis of this lemma is satisfied. Doing this, we replace each time $k$ with $lk$. The permutation~$k$ obtained at the end has the property that the $\Delta$-coordinate of~$k$ does not belong to $K^\Delta$, only if $\Delta\in\orb(\ov S)$.\eprf

\thrml{230419b}
In the  notation and assumption of Theorem~{\em\ref{300419o}}, suppose that the section~$\ov S$ is nonsolvable. Then
\nmrt
\tm{i}  $S^\Delta$ is nonregular for all $\Delta\in\orb(S)$,
\tm{ii} the standard equivalence relations for $S$ and $\ov S$ coincide,
\tm{iii} $\ov S=\soc(\ov S)S$.
\enmrt
\ethrm
\proof To prove (i), we note  that for each $\Delta\in\orb(S)$,  the group $K^\Delta=\AGL_1(\Delta)$ is the normalizer of $C(\Delta)\le S^\Delta$ (Theorems~\ref{181219a} and~\ref{300419o}). It follows that there is $\ov k\in \ov S$ such that $\ov k\ovi{\Delta}$ lies in $K^\Delta$ but not in $C(\Delta)$ (the latter holds, because $p\ge5$). By  Corollary~\ref{060519c}, $\ov k$ can be lifted to $k\in \ov K$ such that $k_\Delta\in K^\Delta$ for all $\Delta$-coordinates of~$k$. It follows that $k\in W_F(K)$.  By Lemma~\ref{251019a}(ii), this implies that $k\in K$ and, therefore, $\ov k\in S$. The group $S^\Delta$ is nonregular, because $\ov k\ovi{\Delta}\in S^\Delta\setminus C(\Delta)$. \medskip

To prove (ii), we apply Lemma~\ref{010418b} for $\Omega=\Omega_S$, $e=e^{}_S$, and
$$
(G,L)\in\{(K^{\Omega_S},S),\ (\ov K\ovi{\Omega_S},\ov S)\}.
$$
It is easily seen that conditions (i) and (ii) of this lemma are satisfied for both pairs. Next, any $\Delta\in\orb(S)$ is of prime cardinality (Theorem~\ref{181219a}). Therefore, the group $S^\Delta$ is primitive. It is nonregular by statement~(i). Since $S^\Delta\le\ov{S}\ovi{\Delta}$,  the group $\ov{S}\ovi{\Delta}$ is primitive and nonregular.  Consequently,
condition~(iii) of  Lemma~\ref{010418b} is also satisfied. Thus, statement~(ii)  follows from Lemma~\ref{010418b}.\medskip

To prove (iii), it suffices to check that $\soc(\ov S)S\ge \ov S$. Let $\ov k\in\ov S$, and $k\in\ov K$ be such that $k^{\Omega_S}=\ov k$ and $k_\Delta\in K^\Delta$ for all $\Delta\not\in\orb(\ov S)$ (see Corollary~\ref{060519c}). By virtue of Lemma~\ref{040519c}, there is $l\in\ov K$ such that $l^{\Omega_S}\in\soc(\ov S)$ and $(lk)_\Delta\in K^\Delta$ for all $\Delta\in\orb(\ov S)$ (if $k_\Delta\in K^\Delta$ for all $\Delta\in\orb(\ov S)$, then one can take $l=1$). It follows that $lk\in \ov K\cap W_F(K)=K$ (Lemma~\ref{251019a}(ii)). Thus,
$$
\ov k=(l^{-1}lk)^{\Omega_S}=(l^{-1})^{\Omega_S}(lk)^{\Omega_S}\in \soc(\ov S) S,
$$
as required.\eprf\medskip

{\bf Proof of Theorem~\ref{200318a}.} Let $G$ be a supersolvable group. Then there exists a maximal normal $G$-flag~$F_0$ (Lemma~\ref{121219a}). Let $F$ be a maximal normal $K$-flag extending~$F_0$, where $K$ is the relative closure of $G$ with respect to~$F_0$ (Corollary~\ref{181219i}). By Lemma~\ref{170419b},  $F$ is also a normal $\ov K$-flag. Therefore, the series
$$
\{\id_\Omega\}=\ov K_{e_0}<\ov K_{e_1}<\cdots<\ov K_{e_m}=\ov K
$$
is normal, where $m$ is the length of~$F$. It follows that all composition factors of $\ov K$ are among the composition factors of the $F$-sections of~$\ov K$.  Let $\ov S$ be an $F$-section of~$\ov K$. Then from Theorem~\ref{230419b}(iii), it follows that all nonabelian composition factors of $\ov S$ are among those of $\soc(\ov S)$. By Theorem~\ref{300419o}(i)(b), each of them is an alternating group of prime degree. Since $\ov G=\ov K$ (Lemma~\ref{251019a}(i)), we are done.\eprf\medskip

We complete the subsection by one more corollary of the reduction lemma, to be used for verifying  the correctness of the Main Algorithm in Subsection~\ref{010519d}.

\thrml{200419d}
Assume that  for every $F$-section $S$ of $K$ such that $\ov S$ is nonsolvable, we are given a set $X_S\subseteq \ov K_{e_S}$ with $\grp{K_{e_S},X_S}^{\Omega_S}=(\ov K_{e_S})\ovi{\Omega_S}$. Then
$$
\ov K=\grp{K,X},
$$
where $X$ is the union of all $X_S$.
\ethrm
\proof It suffices to verify that $\ov K\le \grp{K,X}$. Assume on the contrary that 
\qtnl{131219a}
k\in \ov K\qaq k\not\in  \grp{K,X}.
\eqtn  
If $k_\Delta\in K^\Delta$ for all~$\Delta$, then $k\in W_F(K)$. By Lemma~\ref{251019a}(ii), this implies that $k\in \ov K\cap W_F(K)=K$,
a contradiction.\medskip 

Now we may assume that  for $k$, the hypothesis of Lemma~\ref{040519c} is satisfied for some $F$-section $\ov S$   and $\Delta\in\orb(\ov S)$, and also  the number $i(k)=i(S)$ is maximal possible. Without loss of  generality, we may also assume that $i(k)$ is minimal possible among the~$k$ satisfying~\eqref{131219a}.\medskip 

Let $l\in \ov K$ be as in the conclusion of Lemma~\ref{040519c}. Then $l\in \ov K_{e_S}$. By the theorem hypothesis, this implies that there exists $l'\in\grp{K_{e_S},X_S}$ such that
$$
l^{\Omega/e_{i-1}}=(l')^{\Omega/e_{i-1}}.
$$
Set $k'=l'k$. Then  clearly, $k\in\grp{K,X}$ if and only if $k'\in\grp{K,X}$. Observe that~$k'$  belongs to $\ov K$ and acts on $\Omega/e_S$ as~$k$ does. Moreover, for each $\Delta\in\orb(\ov S)$,
$$
k'_\Delta=(l'k)_\Delta=(l')_\Delta k_\Delta=l_\Delta k_\Delta=(lk)_\Delta,
$$
where the second equality holds because $(l')^{\Omega/e_S}=1$ and one can apply formula~\eqref{040519a}. Thus, $k'_\Delta=(lk)_\Delta\in K^\Delta$ by Lemma~\ref{040519c}. It follows that $i(k')<i(k)$, which contradicts the choice of~$k$.\eprf

\subsection{The generating sets $X_S$.}\label{030519k}  Let $S$ be an $F$-section of~$K$. By Theorems~\ref{181219a} and~\ref{300419o}, the group $S$ is plain and each transitive constituent of~$S$ is a primitive group of prime cardinality.\medskip

In what follows, we assume that each transitive constituent of~$S$ is nonregular (by Theorem~\ref{230419b}(i), this condition is always satisfied  if the section $\ov S$ is nonsolvable). Consequently, a plain bijection $f:\Delta\to\Gamma$ is unique for all $\Delta,\Gamma\in\orb(S)$ such that $\Delta\sim\Gamma$ (Lemma~\ref{240318a}(ii)).\medskip

Given a class $\Lambda$ of the standard equivalence relation for~$S$, we define a permutation $x=x_\Lambda\in\sym(\Omega_S)$ by the action on the $S$-orbits $\Delta$:
\nmrt
\tm{C1} $x^\Delta$ is a $3$-cycle if  $\Delta\in \Lambda$,
\tm{C2} $x^\Gamma=f^{-1}x^\Delta f$ if $\Delta,\Gamma\in\Lambda$, where $f:\Delta\to\Gamma$ is the plain bijection,
\tm{C3} $x^\Delta=\id_\Delta$ if $\Delta\not\in\Lambda$.
\enmrt
A role playing by these very special permutations is clarified in the following statement.

\lmml{230419w}
In the above notation, assume that $\ov S$ is nonsolvable. Let $\ov X_S$ be the set of all permutations~$x_\Lambda$, where $\Lambda$ runs over the classes  of the standard equivalence relation for the group~$S$. Then
\qtnl{010519w}
\grp{S,\ov X_S}=\ov S.
\eqtn
\elmm
\proof From Theorems~\ref{300419o}(i)(b) and~\ref{230419b}(ii), it follows that $\ov X_S\subseteq\ov S$. Therefore, $\grp{S,\ov X_S}\le\ov S$. To prove the reverse inclusion, it suffices to verify that $\soc(\ov S)\le \grp{S,\ov X_S}$ (Theorem~\ref{230419b}(iii)).  However, given $\Delta\in\orb(S)$ the group $S^\Delta$ is primitive, the set $(\ov X_S)^\Delta$ contains a $3$-cycle $x^\Delta$, and $\soc(\ov S)^\Delta=\alt(\Delta)$ (Theorems~\ref{181219a} and~\ref{300419o}). By \cite[Theorem~13.3]{Wie1964}, this implies that 
$$
\soc(\ov S)^\Delta\le \grp{S,\ov X_S}^\Delta.
$$
Furthermore, $\soc(\ov S)$ is a plain group  and  the standard equivalence relation for $\soc(\ov S)$ coincides with that for $\ov S$ (Theorem~\ref{300419o}(ii)) and hence for~$S$ (Theorem~\ref{230419b}(ii)). This completes the proof  by statements~(i) and~(iii) of Lemma~\ref{240318a}.\eprf\medskip

Given an element $x\in \ov X_S$, we define a permutation $k=k_x$ on $\Omega$,  the standard representation  (see Subsection~\ref{131219o}) of  which is uniquely determined by the following three conditions:
\nmrt
\tm{K1} $k^{\Omega/e_S}=\id_{\Omega/e_S}$,
\tm{K2} $k^{\Omega_S}=x$,
\tm{K3} $k_\Delta=\id_\Delta$, $\Delta\in\orb(T)$, where $T$ is an $F$-section of~$K$ with $i(T)< i(S)$.
\enmrt
From Lemma~\ref{230419w}, it follows that $k\in W_F(\ov K)$.  

\thrml{230419x}
In the above notation, assume that $\ov S$ is nonsolvable and $x\in \ov X_S$. Let
$$
C_x:=W_ek_x,
$$
where $W=W_F(K)$ and $e=e_{i(S)-1}$. Then the intersection $\ov K\cap C_x$ is not empty and $(\ov K\cap C_x)^{\Omega_S}=\{x\}$.
\ethrm
\proof  Since $\ov X_S\subseteq\ov S$ (see~\eqref{010519w}), there is   $k\in\ov K_{e_S}$ such that $k^{\Omega_S}=x$. By Corollary~\ref{060519c}, we may assume that $k_\Delta\in K^\Delta$ for all $\Delta\in\orb(T)$, where $T$ is an $F$-section of~$K$ with $i(T)< i(S)$. But then
$$
kk_x^{-1}\in W_e,
$$
which proves that $\ov K\cap C_x$ is not empty. Since $(W_e)^{\Omega_S}=\{\id_S\}$, the second statement immediately follows from the first one.\eprf\medskip

An $F$-section $S$ of the group~$K$ is said to be {\it feasible} if  each transitive constituent of $S$ is a nonregular group of prime degree at least~$5$. A set  $X_S$ is  called a {\it nonsolvability certificate} for~$S$ if exactly one of the following two conditions is satisfied:
\nmrt
\tm{i} $X_S=\varnothing$ if $S$ is not feasible or there exists $x\in\ov X_S$ such that $\ov K\cap C_x=\varnothing$, 
\tm{ii} $X_S$ is  a full set of distinct representatives of the family $\{\ov K\cap C_x:\ x\in\ov X_S\}$,
\enmrt
where  $\ov X_S$ and $C_x$ are as in Lemma~{\rm\ref{230419w}} and~Theorem~{\rm\ref{230419x}}, respectively.

\crllrl{270419a}
Let  $S$ be an $F$-section of $K$. Then  $\ov S$ is nonsolvable if and only if $X_S\ne\varnothing$.
\ecrllr
\proof The necessity follows from Theorem~\ref{300419o}(i)(b) and Theorem~\ref{230419x}. To prove the sufficiency, let $x\in\ov X_S$ (note that  $\ov X_S$ is well defined, because $S$ is feasible). Then  there exists $k\in \ov K$ such that $k^{\Omega_S}=x$; in particular,  $k^\Delta\in\ov S\ovi{\Delta}$ is a $3$-cycle for some $\Delta\in\orb(S)$. On the other hand, $|\Delta|\ge 5$ and hence the group $\AGL_1(\Delta)$ contains no $3$-cycles. Thus by Theorem~\ref{300419o}(i)(a), the group $\ov S\ovi{\Delta}$, and hence $\ov S$, is nonsolvable.\eprf

\crllrl{210419f}
Let  $S$ be an $F$-section of $K$ such that  $\ov S$ is nonsolvable.  Then 
$$
X_S\subseteq \ov K_{e_S}\qaq \grp{K_{e_S},X_S}^{\Omega_S}=(\ov K_{e_S})^{\Omega_S}.
$$
\ecrllr
\proof The first inclusion immediately follows from the definition of~$X_S$. Further, by  Theorem~\ref{230419x}  the intersection $\ov K\cap C_x$ is not empty for all $x\in \ov X_S$. It follows that $(X_S)^{\Omega_S}=\ov X_S$. Since $S=(K_{e_S})\ovi{\Omega_S}$ and $\ov S=(\ov K_{e_S})\ovi{\Omega_S}$, the required equality is a consequence of Lemma~\ref{230419w}.~\eprf

\section{Proof of Theorem~\ref{200318}}\label{240318u}

The Main Algorithm and its analysis is given at the end of the section. We start with description of auxiliary procedures used  in the Main Algorithm.

\subsection{Finding a maximal normal flag.}\label{010519a} It is known (see \cite[p.~49]{Seress}) that a minimal normal subgroup of a permutation group of degree~$n$ can be constructed in time $\poly(n)$. It follows that series~\eqref{290319a} can also be found efficiently. This proves the statement below.

\lmml{290319b}
A maximal normal flag of a permutation group of degree $n$ can be constructed in time~$\poly(n)$.
\elmm

\subsection{Finding the relative closure.} Let $G\le\sym(\Omega)$ and $F$ a $G$-flag.  To find the majorant $W_F(G)$, one needs first to identify each $G$-orbit $\Delta$ with the Cartesian product of quotient sets defined by the flag~$F_\Delta$ (see formula~\eqref{270419i}), and to construct a generator set $X_\Delta$ of the iterated wreath product~\eqref{27042019j} for $G=G^\Delta$ and suitable $\Delta_i$'s. This can be done efficiently by  standard permutation group algorithms. Then the generator set of $W_F(G)$ is obtained by taking the  union of all the $X_\Delta$.  \medskip

The  relative closure $\ov G\cap W_F(G)$ can be interpreted as the group 
$\aut(\fG)\cap W_F(G)$, where $\fG$ is the arc-colored graph with vertex set $\Omega$, the color classes of which are exactly the 2-orbits of~$G$.  Now if $G$ is a solvable group, then the group $W_F(G)$ is also solvable (Lemma~\ref{251019a}(iii)). Therefore, a generator set for $\ov G\cap W_F(G)$ can be constructed in time $\poly(n)$ by the Babai--Luks algorithm~\cite[Corollary~3.6]{BL83}.

\lmml{210419e}
The relative closure (with respect to a flag) of a solvable permutation group of degree $n$ can be constructed in time $\poly(n)$.
\elmm

\subsection{Extending of a flag.} The following algorithm provides  a constructive version of Corollary~\ref{181219i}. Below, given a flag and its element $e=e_{i-1}$, we set $\ov e=e_i$.

\medskip
\centerline{\bf Extending to maximal normal flag}
\medskip

\noindent {\bf Input:} a  group $K\le\sym(\Omega)$ and a normal $K$-flag $F$. 

\noindent {\bf Output:} a maximal normal $K$-flag $F'$ extending $F$.

\medskip

\noindent{\bf Step 1.} Set $F':=F$.\smallskip

\noindent{\bf Step 2.} While there is $e\in F'$ and a non-singleton $\Lambda\in\orb(K\ovi{\Omega/e})$  such that $((K_{\ov e})\ovi{\Omega/e})_\Lambda$ contains a proper minimal normal subgroup $N$ of $K\ovi{\Omega/e}$, set 
$$
F':=F'\,\cup\,\{e'\},
$$ 
where $e'$ is the equivalence relation on~$\Omega$, the classes of which on $\Omega/e$ are the $N$-orbits.

\noindent{\bf Step 3.} Output $F'$.
\eprf\medskip

\thrml{270419o6}
Assume that $K$ is the relative closure of a supersolvable group~$G$ with respect to~$F$, and  $F$ is a maximal normal~$G$-flag. Then the above algorithm correctly finds a maximal normal $K$-flag $F'$, extending~$F$,  in time $\poly(n)$.
\ethrm
\proof Each iteration at Step~2 increases the length of the flag~$F'$. Therefore, the number of iterations is at most~$n$. The cost of the iteration consists of two parts: finding the pointwise stabilizer$((K_{\ov e})\ovi{\Omega/e})_\Lambda$  and  constructing minimal normal subgroup~$N$.  Thus the running time of the algorithm is $\poly(n)$, see~\cite{Seress}. Finally, the $K$-flag~$F'$ defined at Step~2 is normal by the construction, and the resulted flag is maximal normal by Corollary~\ref{181219i}.\eprf

\subsection{An nonsolvability certificate.} In this subsection, we present the following algorithm finding  the nonsolvability certificate for a given section of a permutation group with respect to a given flag.

\medskip
\centerline{\bf Finding nonsolvability certificate}
\medskip

\noindent {\bf Input:} a  group $K\le\sym(\Omega)$, a normal $K$-flag $F$,  and a  $F$-section $S$ of~$K$.

\noindent {\bf Output:} the nonsolvability certificate $X_S$.

\medskip

\noindent{\bf Step 1.} If $S$ is not feasible, then output $\varnothing$.\smallskip

\noindent{\bf Step 2.} If $S$ is not plain, then  output $\varnothing$; else construct the standard equivalence relation of~$S$.

\noindent{\bf Step 3.} Construct the set $\ov X_S$ defined in Lemma~\ref{230419w}.\smallskip

\noindent{\bf Step 4.} Set $X_S:=\varnothing$. For each $x\in\ov X_S$,

{\bf Step 4.1.} Find the intersection $Y:=\ov K\cap C_x$.

{\bf Step 4.2.} If $Y=\varnothing$, then output $\varnothing$.

{\bf Step 4.3.} Add to $X_S$ one (arbitrary) element of~$Y$.\smallskip

\noindent{\bf Step 5.} Output $X_S$.
\eprf\medskip

\thrml{270419o}
Assume that $K$ is the relative closure of a supersolvable group, and  $F$ is a maximal normal~$K$-flag. Then the above algorithm correctly finds the nonsolvability certificate in time $\poly(n)$.
\ethrm
\proof The correctness of the algorithm follows from the definition of the nonsolvability certificate. To estimate the running time, we note that Steps~1 and~2 require  inspecting, respectively, the orbits and $2$-orbits of~$S$, and this can be done by standard permutation group algorithms. The complexity of Step~3 is bounded by a polynomial in
$$
|\ov X_S|=|\orb(S)/\sim|\le n,
$$
where $\sim$ is the standard equivalence relation for~$S$ (recall that $S$ is a plain group by Theorem~\ref{300419o}(ii)).\medskip

At Step~4, one needs to find $|X_S|$ times the intersection $Y$ of two cosets $\ov K$ and~$C_x$. Note that the latter is a coset of a solvable group contained in~$W_F(K)$ (Theorem~\ref{230419x}). Therefore, finding the intersection can be interpreted as finding the automorphisms of the arc-colored graph associated with~$K$, and hence can be implemented in time $\poly(n)$ by the Babai--Luks algorithm~\cite{BL83}. Thus Step~4 and  the whole algorithm run in polynomial time.\eprf

\subsection{The Main Algorithm.}\label{010519d} In this subsection, we complete the proof of Theorem~\ref{200318}.

\medskip
\centerline{\bf Main Algorithm}
\medskip

\noindent {\bf Input:} a supersolvable group $G\le\sym(\Omega)$.

\noindent {\bf Output:} the $2$-closure $\ov G$.

\medskip

\noindent{\bf Step 1.} Find a maximal normal flag $F$ of $G$ (Lemma~\ref{290319b}).\smallskip

\noindent{\bf Step 2.} Find the relative closure $K$ of $G$ with respect to~$F$ (Lemma~\ref{210419e}).\smallskip

\noindent{\bf Step 3.} Find a maximal normal $K$-flag $F'$ extending~$F$ (Theorem~\ref{270419o6}).\smallskip

\noindent{\bf Step 4.} Set $X:=\{\id_\Omega\}$. For each $F'$-section $S$ of $K$, find the nonsolvability certificate~$X_S$, and put $X:=X\,\cup X_S$ (Theorem~\ref{270419o}).\smallskip

\noindent{\bf Step 5.} Output $\ov G:=\grp{K,X}$.
\eprf\medskip

\thrml{140318a}
The Main Algorithm correctly finds the group $\ov G$ in time~$\poly(n)$.
\ethrm
\proof The fact that the Main Algorithm is polynomial-time follows from Lemmas~\ref{290319b} and~\ref{210419e}, and Theorems~\ref{270419o} and~\ref{270419o6}. Next, let $S$ be a section appearing at Step~4 and such that $X_S\ne\varnothing$. Then the section $\ov S$ is nonsolvable (Corollary~\ref{270419a}) and hence the conditions of Theorem~\ref{200419d} are satisfied (Corollary~\ref{210419f}). By this theorem,  the output group coincides with $\ov K$. Since also $\ov K=\ov G$ (Lemma~\ref{251019a}(i)), the output  of the Main Algorithm is correct.\eprf

\end{document}